\documentclass[12pt]{amsart}

\usepackage{amssymb}

\usepackage{enumitem}

\usepackage{graphicx}

\makeatletter
\@namedef{subjclassname@2020}{%
  \textup{2020} Mathematics Subject Classification}
\makeatother

\usepackage[T1]{fontenc}

\newtheorem{theorem}{Theorem}[section]

\newtheorem{claim}[theorem]{Claim}

\theoremstyle{definition}

\newtheorem{remark}[theorem]{Remark}

\numberwithin{equation}{section}

\begin{document}


\baselineskip=17pt


\title[On product of two non-trivial countable-dimensional continua]{No product of two non-trivial countable-dimensional continua maps lightly into any of the factors}
\author{Roman Pol and Miros{\l}awa Re\'{n}ska}

\address{Institute of Mathematics\\ University of Warsaw\\Banacha 2\\ 02-097 Warszawa, Poland}
\email{r.pol@mimuw.edu.pl}

\address{Faculty of Mathematics and Information Science\\Warsaw University of Technology\\ 
Koszykowa 75\\00-662 Warszawa, Poland}
\email{miroslawa.renska@pw.edu.pl}


\begin{abstract}
We shall prove that if $X$, $Y$ are compact metrizable spaces of positive dimension 
and $h:X\times Y \to X$ is a continuous map with zero-dimensional fibers then $X$ 
contains a non-trivial continuum without one-dimensional subsets; in particular $X$ is 
not a countable union of zero-dimensional sets, which 
provides a negative answer to a question of 
J. Dud\'ak and B. Vejnar \cite{DV}
\end{abstract}

\maketitle
\section{Introduction}

All our spaces are metrizable separable and a compactum means a compact space. Our terminology follows  \cite{En1} and \cite{En2} ,
$I=[0,1]$ is the unit interval and $I^{\mathbb N}$ is the Hilbert cube.

Countable-dimensional spaces are the spaces which are countable unions of zero-dimensional sets.

A mapping between compacta is light if it is continuous and no fiber of the map contains a non-trivial 
continuum (i.e., the fibers are zero-dimensional).

Since each countable-dimensional compactum of positive dimension contains a one-dimensional compactum, cf. 
\cite{HW}, Ch. IV, \textsection 6, the following result provides a negative answer to a question asked by 
J. Dud\'ak and B. Vejnar \cite{DV}, Question 3.8.

\medskip

\begin{theorem} Let $X$, $Y$ be compacta of positive dimension such that there is a light map 
$h:X\times Y \to X$. Then $X$ contains a compactum of positive dimension without one-dimensional 
subsets. 
\end{theorem}

\medskip

Non-trivial continua without one-dimensional subsets (equivalently - without subsets of positive finite dimension) 
were first constructed by \cite{Wa}, cf. \cite{En2}, sec. 5.2.

If $h$ in Theorem 1.1 is an embedding, the compactum $X$ is strongly infinite-dimensional, i.e., 
$X$ admits an essential map onto $I^{\mathbb N}$, cf. Remark 3.3 and [He].

However, J. van Mill \cite{vM} constructed a Peano continuum $X$ homeomorphic to  $X^2$ 
but not to $X^{\mathbb N}$.

The scheme of our proof is the following. We construct by transfinite induction compacta 
$E_{\alpha }$, $\alpha < \omega_1$, which admit essential maps onto $AR$-compacta 
$H_{\alpha }$ defined by D. W. Henderson \cite{He}, cf. sec. 2 (A), and each $E_{\alpha }$ 
maps lightly into $X$. Then we consider a compactum $Z$ which contains a copy of every 
$E_{\alpha }$ and maps lightly into $X$. The compactum $Z$ maps essentially onto 
each Henderson compactum $H_{\alpha }$, $\alpha < \omega_1$, hence by \cite{Po1}, 
$Z$ is strongly infinite-dimensional and its light image satisfies the assertion of Theorem 1.1.

We shall prove Theorem 1.1 in sec. 3, after explaining in sec. 2 some background needed 
for the proof.

In the last section, we give some more information related to Theorem 1.1, in the case 
the map $h$ is an embedding.

\section{Smirnov's compacta, Henderson's compacta and essential maps}. 

{\bf (A)} Smirnov's compacta $S_{\alpha}$, $\alpha < \omega_1$, are defined by 
transfinite induction, cf. \cite{En1}, 7.1.33, \cite{Po1},  \textsection 2, sec. 2: 
$S_1 = I$, $S_{\alpha + 1} = S_{\alpha} \times I$ and for any limit $\alpha < \omega_1$, 
 $S_{\alpha}$ is the one-point compactification of the free union  $\bigoplus\limits_{\beta < \alpha} S_{\beta}$.

We shall identify  $S_{\alpha}$ with the subspace  $S_{\alpha} \times \{ 0 \}$ of  $S_{\alpha + 1}$, and for 
any limit $\alpha$ and $\beta < \alpha$,  $S_{\beta}$ is embedded in  $S_{\alpha}$ as an open subspace.

Open components of  $S_{\alpha}$ are closed $n$-cells, i.e., homeomorphic copies of $n$-cubes $I^n$. 

In the process of construction of Smirnov's compacta we shall also fix, for each open component 
$C$ of  $S_{\alpha}$, a homeomorphism $e_C:C \to I^n$, keeping an agreement that, for any open component 
$D = C \times I$ of  $S_{\alpha + 1}$, $C$ being an open component of  $S_{\alpha}$, 
$e_D = e_C \times id_I : D \to I^{n+1}$ ($ id_I$ is the identity on $I$).

If $\alpha$ is a limit ordinal and $C$ is an open component of  $S_{\alpha}$, 
$C$ is also an open component of  the open copy $S_{\beta}$ in $S_{\alpha}$ containing $C$, 
$\beta < \alpha$, and we let $e_C$ be the map already fixed for this component of $S_{\beta}$.

For any open component $C$ of $S_{\alpha}$, $\partial C = e_C^{-1} (\partial I^n)$ is the geometric 
boundary of $C$.

D. W. Henderson  \cite{He}, cf. 
 \cite{Po1},  \textsection 3, sec. 1, extended by transfinite induction each $S_{\alpha}$ to an $AR$-compactum 
$H_{\alpha} \supset S_{\alpha}$ such that for each open component $C$ of $S_{\alpha}$, $C\setminus \partial C$ 
is open in $H_{\alpha}$, and distinguished the closed set 
$\partial  H_{\alpha} = H_{\alpha} \setminus \bigcup \{ C\setminus \partial C : C\ {\rm is\ an\ open\ component\ of }\  S_{\alpha} \}$.

\medskip

{\bf (B)} Let us recall that a continuous map $f:K \to C$ of a compactum onto a closed cell is essential, if for any continuous 
$g:K \to C$ coinciding with $f$ on $f^{-1} (\partial C)$, $g(K) = C$, cf. \cite{En2} (in \cite{GT}, 
sec. 5, essential maps onto closed cells are called $AH$-essential).

A continuous map $f:K \to I^{\mathbb N}$ of a compactum is essential if for any projection 
$p_n:I^{\mathbb N} \to I^n$ onto the first $n$ coordinates, $p_n\circ f :K \to I^n$ is essential, 
cf.  \cite{Al},  \cite{He}.

The compacta admitting essential maps onto the Hilbert cube are exactly the strongly 
infinite-dimensional ones, cf.  \cite{Al},  \cite{En2}.

Henderson  \cite{He} called a continuous map $f:K\to H_{\alpha}$ of a compactum onto the Henderson 
compactum essential if for each continuous $g:K \to H$ coinciding with $f$ on $f^{-1} (\partial H_{\alpha})$, 
$g(K) = H_{\alpha}$, cf. (A).

\medskip 

\begin{remark}
By  \cite{He}, Proposition 3, a continuous surjection  $f:K\to H_{\alpha}$ is essential if and only if 
for each open component $C$ of the Smirnov compactum  $S_{\alpha} \subset  H_{\alpha}$, the restriction 
$f | {f^{-1}(C)}:f^{-1}(C) \to C$ is an essential map onto the cell.
\end{remark}

\medskip 

Our reasonings will be based on the following result.

\medskip

\begin{theorem} ( \cite{Po1},  \textsection 3, Theorem 2.1). 
A compactum which admits for every $\alpha < \omega_1$ an essential map onto 
the Henderson compactum $H_{\alpha}$, admits also an essential map onto  
$ I^{\mathbb N}$.
\end{theorem}

\medskip 

{\bf (C)} The following theorem of Holszty\'{n}ski will be an important ingredient of our proof 
( Holszty\'{n}ski calls in \cite{Ho} essential maps onto cubes ''universal'', cf. \cite{GT}, Theorem 5.4).

\medskip 

\begin{theorem} ( \cite{Ho}). For any continuous surjections $u_i:T_i \to I$, $i\leq n$, 
$T_i$ being continua, the product map 
$u_1\times u_2 \times \ldots \times u_n: \prod\limits_{i=1}^{n} T_i \to I^n$ is essential.
\end{theorem}

\section{Proof of theorem 1.1.}
{\bf (A)} We adopt the notation from sec. 2 (A), We shall say that a continuous map 
$f:K \to S_{\alpha}$ of a compactum onto the Smirnov compactum is determined by products 
if for every open component $C$ of $S_{\alpha}$ there is a homeomorphism 
$g:f^{-1} (C) \to T_1 \times  T_2 \times \ldots  \times  T_n$, $n= {\rm dim } C$, onto 
the product of continua  $ T_i$, and continuous surjections $u_i:T_i \to I$,  $u_i^{-1} (0)$ 
being a singleton, such that 
$f | {f^{-1}(C)} = e_C^{-1} \circ (u_1\times u_2 \times \ldots \times u_n ) \circ g$.

Let $\alpha = \beta + 1$, $\beta \geq 1$, and let us recall that $S_{\beta }$ is identified 
with  $S_{\beta } \times \{ 0 \} \subset S_{\alpha}$. Then, whenever $f:K \to S_{\alpha}$ 
 is determined by products, so is the restriction 
$f | f^{-1}(S_{\beta }):  f^{-1}(S_{\beta }) \to  S_{\beta }$.

Indeed, if $D$ is an open component of $ S_{\beta }$, $C= D \times I$ is an open component 
of $S_{\alpha}$  and in the above notation, $g$ maps $f^{-1} (D)$ onto 
$T_1 \times  T_2 \times \ldots  \times  T_{n-1} \times \{ a\}$, $ \{ a\}= u_n^{-1}(0)$.

Clearly, if $f:K \to S_{\alpha}$  is determined by products, $\alpha$ is limit, $\beta < \alpha$, and 
$S_{\beta }$ is the open copy of the  Smirnov compactum in  $S_{\alpha}$, then 
$f | f^{-1}(S_{\beta }):  f^{-1}(S_{\beta }) \to  S_{\beta }$ is also determined by products.

This implies readily by transfinite induction on $\alpha$ that if a compactum $K$ admits a 
continuous map onto  $S_{\alpha}$  determined by products and $\beta < \alpha$, 
then there is a compactum $L$ in $K$ admitting a continuous map onto  
$S_{\beta}$  determined by products. 

Finally, by the  Holszty\'{n}ski Theorem 2.3 and Remark 2.1, if 
$f:K \to H_{\alpha}$ is a continuous map of a compactum onto the Henderson compactum 
and for some compactum $L \subset K$, the restriction $f | L : L \to  S_{\alpha}$ is determined by products, 
then $f$ is an essential map onto $H_{\alpha}$.

\medskip

{\bf (B) } Let $X$, $Y$ be compacta of positive dimension and let $h: X \times Y \to X$ be a light map. 
We shall fix metrics $d_X$, $d_Y$ generating the topology of $X$ and $Y$, respectively, and we endow 
$ X \times Y $ with the metric 

\smallskip

$d_{ X \times Y } \bigl( (x_1,y_1),  (x_2,y_2) \bigr) = {\rm max}  \bigl( d_X(x_1,x_2), d_Y (y_1,y_2) \bigr)$.

We shall use the uniform continuity of $h$ with respect to the metric $d_{ X \times Y }$ in 
$ X \times Y$ and the metric $d_ X$ in $X$.

We shall check by transfinite induction the following claim.

\begin{claim}
For each  $\alpha < \omega_1$ and $\varepsilon > 0$, there is a compactum $K$ and a light map 
$f:K \to X$ such that $K$ has a continuous map onto $S_{\alpha }$ determined by products and 
${\rm diam} f(K) < \varepsilon$.
\end{claim}

\medskip

The claim is true for $\alpha = 1$. Indeed, given  $\varepsilon > 0$, we pick a non-trivial continuum 
$K \subset X$ with ${\rm diam} K < \varepsilon$ and $f$ is the identity on $K$ (recall that 
$S_1 = I$).

Assume that Claim 3.1 is true for all $\beta < \alpha$, $\alpha$ being a countable ordinal.

Assume first that $\alpha = \beta + 1$. Let us fix $\varepsilon > 0$. and let $\delta > 0$ 
be such that $h$ takes sets in $X \times Y$ of diameter $< \delta$ to sets in $X$ of 
diameter $< \varepsilon$. By the inductive assumption, there is a compactum $L$, a 
continuous surjection $u:L \to S_{\beta }$ determined by products and a light map 
$g:L \to X$ with ${\rm diam}\bigl( g(L) \bigr)< \delta$.
Let $T$ be a continuum in $Y$ with ${\rm diam}T < \delta$ and let $w:T \to I$ be a continuous 
surjection such that $w^{-1} (0)$ is a singleton.

Then the continuous map $u:L \times T \to  S_{\beta } \times I =  S_{\beta + 1 } =  S_{\alpha }$ 
is determined by products and ${\rm diam} \bigl( h\bigl( g(L) \times T \bigr) \bigr) < \varepsilon$.
The continuous map $f=h\circ (g\times id_T) : L \times T \to X$ is light ($ id_T$ is the identity on $T$), 
as the composition of two light maps and  
${\rm diam} f(L \times T) < \varepsilon$.

Assume now that $\alpha$ is a limit ordinal.
Let $\varepsilon > 0$ and $\delta > 0$ be as in the previous case.
Let us list all ordinals less than $\alpha $ as $\alpha_1, \alpha_2, \ldots $ 
and let us choose ordinals $\beta_1 < \beta_2 < \ldots $ with $\alpha_n \leq \beta_n < \alpha$. 

The inductive assumption provides compacta $L_n$ which have continuous maps onto $ S_{\beta_n }$ 
determined by products, and light maps $f_n:L_n \to X$ with ${\rm diam} f_n(L_n) < {1\over n}$.

By compactness of $X$, passing if necessary to a subsequence, we can ensure that there is $a\in X$ such that 
${\delta \over 2} > {\rm diam} \bigl(  \{ a \} \cup f_n(L_n)   \bigr) \to 0$.

The observation made in (A) allows one to pick for each $n$ a compactum $M_n$ in 
$L_{n }$ which has a continuous map onto $ S_{\alpha_n }$ determined by products.

The one-point compactification $K$ of the free union $\bigoplus\limits_n M_n$ has a continuous map onto 
$S_{\alpha }$ determined by products, $S_{\alpha }$ being the one-point compactification of the free union 
$\bigoplus\limits_n S_{\alpha_n }$ (the map takes $M_n$ onto $S_{\alpha_n }$ and the point 
at infinity of $K$ to the point at infinity of $S_{\alpha }$).

Let $b_n \to b$ be a sequence in $Y$ with all $b, b_1, b_2, \ldots $ distinct and 
${\rm diam} \{ b, b_1, b_2, \ldots \} < \delta$. Then 
$C = \{ (a,b) \} \cup \bigcup\limits_{n=1}^\infty  \bigl(   f_n(M_n)  \times \{ b_n \}  \bigr) $ is 
topologically the one-point compactification of $\bigoplus\limits_n f_n(M_n)$, $(a,b)$ being the point at infinity, 
hence the map $u:K \to C$ taking the point at infinity of $K$ to $(a,b)$ and such that 
$u(x) = \bigl( f_n (x), b_n)$ whenever $x \in M_n$, is a light map. 
The composition $h \circ u : K \to X$ is light 
and ${\rm diam} \bigl( h \circ u  (K) \bigr) < \varepsilon$.

{\bf (C) } Having checked Claim 3.1, we obtain for each $\alpha < \omega_1$ a compactum 
$K_{\alpha } \subset   I^{\mathbb N}$ which has a continuous map onto $S_{\alpha }$ 
determined by products and a light map $f_{\alpha } : K_{\alpha } \to X$ (the additional 
condition involving $\varepsilon > 0$ was only needed to carry out the inductive construction).

We shall 
verify the following claim.

\begin{claim}
There is a compactum $Z$ and a light map 
$w:Z \to X$ such that each compactum $K_{\alpha }$ embeds in $Z$.
\end{claim}

To that end, let us take a continuous surjection 
$\Phi: 2^{\mathbb N} \to \mathcal{K} (I^{\mathbb N} \times X)$  of the 
Cantor set onto the space of non-empty compact sets in $I^{\mathbb N} \times X$, equipped 
with the Vietoris topology, cf. \cite{En1}, 3.12.20, 4.5.23, and let 

$H = \{ (t,x,y): (x,y) \in \Phi (t) \}$, $p(t,x,y) = y$. 

Then 
$H= \bigcup\limits_{t \in 2^{\mathbb N} } \{ t \} \times  \Phi (t) \subset  2^{\mathbb N} \times I^{\mathbb N} \times X $ 
is a compactum and $p:H \to X$ is a continuous map whose restriction to $ \{ t \} \times  \Phi (t)$ is the projection 
of  $\Phi (t)$ onto the second coordinate. 

The graph ${\rm Gr} (f_{\alpha } ) = \{ \bigl( x, f_{\alpha} (x) \bigr) : x \in K_{\alpha } \}$  is a compact 
set in $  I^{\mathbb N} \times X $, hence ${\rm Gr} (f_{\alpha }) =  \Phi (t_{\alpha } ) $ for some 
$t_{\alpha } \in 2^{\mathbb N}$. The compactum $ \{ t_{\alpha } \} \times {\rm Gr} (f_{\alpha })$ 
is homeomorphic to $K_{\alpha } $ and the restriction of $p$ to this compactum is light, as 
$p\bigl( t_{\alpha } , x, f_{\alpha } (x) \bigr) = f_{\alpha } (x)$.

Let $u:H\to Z$, $w:Z \to X$ be the monotone-light factorization of $p$. i.e., the fibers of $u$ are 
connected and the fibers of $w$ are zero-dimensional, cf. \cite{En1}, 6.2.22.
In particular, the map $u$ is obtained by matching each non-trivial continuum in every fiber $p^{-1}(x)$ 
to a point.

Now, any continuum in $H$ is contained in some $ \{ t \} \times  \Phi (t)$, and if a continuum is contained 
in $ \{ t_{\alpha } \} \times {\rm Gr} (f_{\alpha })$ and taken by $p$ to a point, then the 
continuum is a singleton, the restriction of $p$ to $ \{ t_{\alpha } \} \times {\rm Gr} (f_{\alpha })$ 
being light. It follows that $u$ is injective on $ \{ t_{\alpha } \} \times {\rm Gr} (f_{\alpha })$ and in 
effect, $u \bigl( \{ t_{\alpha } \}  \times {\rm Gr} (f_{\alpha }) \bigr)$ is a copy of $K_{\alpha }$ in $Z$.  

\medskip

{\bf (D)} With Claim 3.2 in hand, we are ready to complete the proof of Theorem 1.1. 
Let $w:Z\to X$ be as in Claim 3.2. For each $\alpha < \omega_1$ we have a compactum 
 $L_{\alpha } \subset Z$ and a continuous map  $g_{\alpha } : L_{\alpha } \to S_{\alpha }$ onto the 
Smirnov compactum determined by products.

Let $Z^* = Z \bigoplus 2^{\mathbb N} $ and let  $v_{\alpha } :2^{\mathbb N} \to \partial  H_{\alpha }$ 
be a continuous surjection, $\alpha < \omega_1$. Since $H_{\alpha }$ is an $AR$, the combination of 
 $g_{\alpha }$ and  $v_{\alpha }$ can be extended from $L_{\alpha } \cup 2^{\mathbb N}$ to 
a continuous surjection $\tilde{g}_{\alpha } : Z^* \to H_{\alpha }$ . By the observation at the end of (A), 
 $\tilde{g}_{\alpha }$ is essential. By Theorem 2.2, $Z^*$ is strongly infinite-dimensional, and so is $Z$. 

By \cite{Po2}, 2.3, there is a countable-dimensional $A\subset X$ such that $X\setminus A$ contains no one-dimensional subsets. 
The map $w$ being light, $w^{-1}(A)$ is countable-dimensional and therefore, by \cite{En2}, 6.1.17, there is a 
non-trivial continuum $C\subset Z \setminus w^{-1}(A)$. The non-trivial continuum $w(C)$ contains no one-dimensional 
subsets.

\medskip

\begin{remark} 
If the map $h:X \times Y \to X$ is an embedding, the compacta $K$ constructed in the proof of Claim 3.1 are 
actually subsets of $X$, and hence the reasoning shows that $X$ is strongly infinite-dimensional, cf. Theorem 2.2, 
\cite{Po1},  \textsection 3, Theorem 2.1.
In this case, the space $Z$ constructed in (C) and used in (D) is not needed.
\end{remark}

\section{Comments.}

The proofs of the results stated in this section will be presented by the second author in a separate paper \cite{Re}.

\medskip

{\bf 4.1. Embedding of countable products.} 
If $h:X \times Y \to X$ in Theorem 1.1 is an embedding and for some non-trivial continuum $K$, $Y$ contains
copies of $K$ with arbitrarily small diameters, then 
$K^{\mathbb N}$ embeds in $X$.

\medskip

{\bf 4.2. On the transfinite dimension ind.}
Let $h:X \times Y \to X$ in Theorem 1.1 be an embedding.
Combining the approach of this paper with a subtle theorem of Levin \cite{Le}, one can 
prove that $X$ contains countable-dimensional compacta with arbitrarily large 
transfinite dimension ind, cf. \cite{En2}.

\medskip

{\bf 4.3. On non-compact spaces.}
Let $X$ be a completely metrizable separable space and let 
$h:G\to X$ be a continuous injective map of an open set $G$ in $X \times X$ such that the image 
$h(G)$ is dense in $X$.

If $K$ is a non-trivial continuum such that each non-empty open set in $X$ contains a topological copy of $K$, then 
$K^{\mathbb N}$ embeds in $X$.

\end{document}